\documentclass[11pt]{amsart}
\usepackage{amssymb}
\theoremstyle{plain}
\newtheorem{theorem}{Theorem}
\newtheorem{corollary}{Corollary}

\newtheorem{lemma}{Lemma}
\newtheorem{proposition}{Proposition}
\theoremstyle{example}
\newtheorem{example}{Example}
\theoremstyle{definition}

\theoremstyle{remark}

\numberwithin{equation}{section}
\def\K#1#2{\displaystyle{\mathop K\limits_{#1}^{#2}}}

\setbox0=\hbox{$+$}
\newdimen\plusheight
\plusheight=\ht0
\def\+{\;\lower\plusheight\hbox{$+$}\;}

\setbox0=\hbox{$-$}
\newdimen\minusheight
\minusheight=\ht0
\def\-{\;\lower\minusheight\hbox{$-$}\;}

\setbox0=\hbox{$\cdots$}
\newdimen\cdotsheight
\cdotsheight=\plusheight
\def\cds{\lower\cdotsheight\hbox{$\cdots$}}
\begin{document}
\title[Continued Fractions with Multiple Limits ]
       {Continued Fractions with Multiple Limits}
\author{Douglas Bowman}
\address{ Northern Illinois University\\
   Mathematical Sciences\\
   DeKalb, IL 60115-2888
} 
\email{bowman@math.niu.edu}
\author{J. Mc Laughlin}
\address{Mathematics Department\\
 Trinity College\\
300 Summit Street, Hartford, CT 06106-3100}
\email{james.mclaughlin@trincoll.edu}
\date{March, 31, 2004}
\thanks{The first author's research was partially supported by NSF grant
DMS-0300126.}
\thanks{The second author's research was partially supported by a Trjitzinsky
Fellowship.}
\begin{abstract}
For  integers $m \geq 2$, we study divergent continued fractions
whose numerators and denominators in each of the $m$ arithmetic progressions
modulo $m$ converge. Special cases give, among other things, an infinite
sequence of divergence theorems, the first of which is the classical Stern-Stolz
theorem.

We give a theorem on a general class of Poincar{\'e} type recurrences 
which shows that they
tend to limits when the limits are taken in residue classes and the roots of their
characteristic polynomials are distinct roots of unity.

We also generalize a curious $q$-continued fraction of Ramanujan's
with three limits to a continued fraction with $k$ distinct limit points, $k\geq 2$. 
The $k$ limits are evaluated in
terms of ratios of certain unusual $q$ series.

Finally, we show how to use Daniel Bernoulli's continued fraction in an elementary way
to create analytic  continued fractions with
$m$ limit points, for any positive integer $m \geq 2$.
\end{abstract}

\maketitle

%           Introduction

\section{Introduction } \label{S:intro}

When one studies an infinite process, and it is found not to tend
to a definite limit, an initial instinct is for one to discard the
process as unsuitable. However, there are  cases in which the
divergence occurs in such a controlled way that the process still
retains its utility. Summability provides one example. Another way
in which a divergent process can be useful is if it tends to a
finite number of definite limits for ``nice'' subsequences of its
approximants. This occurs in a natural way in the context of
continued fractions, recurrence sequences and infinite products of
matrices. Here we will make an intensive study of this behavior
when the subsequences are arithmetic progressions of residue
classes modulo $m$.

We begin by reviewing
notation for continued fractions. The symbol $\K{}{}$ is used for continued fractions
in the same way that $\sum$ and $\prod$ are used for series and products,
respectively. Thus,
{\allowdisplaybreaks
\begin{align*}
\K{n = 1}{N}
\displaystyle{\frac{a_{n}}{b_{n}}} \,&:=
\cfrac{a_{1}}{b_{1} + \cfrac{a_{2}}{b_{2} + \cfrac{a_{3}}{b_{3} +
\ldots +  \cfrac{a_{N}^{}}{b_{N} }}}} \\
\phantom{asd} &=  \frac{a_{1}}{b_{1}}
\+
 \frac{a_{2}}{b_{2}}
\+
\frac{a_{3}}{b_{3}}
\+
\cds
\+
\frac{a_{N}}{b_{N}}.
\end{align*}
}
We write $\displaystyle{P_{N}/Q_{N}}$ for the above finite
continued fraction written as a rational function of the variables
$a_{1},...,a_{N},b_{1},...,b_{N}$. $P_{N}$ is the \emph{$N$-th
canonical numerator}, $Q_{N}$ is the \emph{$N$-th canonical
denominator}
and the ratio $P_{N}/Q_{N}$ is the \emph{ $N$-th approximant}.
Let $\hat{\mathbb{C}}$ denote the extended complex plane.

By $K_{n = 1}^{\infty}\displaystyle{a_{n}/b_{n}}$
 we  mean the limit
of the sequence $\{$$\displaystyle{P_{n}/Q_{n}}$$\}$ as \,$n$\,goes to
infinity, if the limit exists.

Two of the most interesting examples of continued fractions with more than one limit are
due to Rogers-Ramanujan \cite{Rcp,Rogers},  and Ramanujan \cite{S88}, respectively:
\begin{equation}\label{RR}
 1+
\frac{q}{1}
\+
\frac{q^{2}}{1}
\+
\frac{q^{3}}{1}
 \+
\frac{q^{4}}{1}
\cds,
\end{equation}
and
\begin{equation}\label{R3}
\frac{-1}{1+q}
\+
\frac{-1}{1+q^2}
\+
\frac{-1}{1+q^3}
\+
\cds .
\end{equation}

Now, it is known that (\ref{RR}) converges for $|q|<1$, while for
$|q|>1$, it tends to two different limits, depending on whether
one considers the sequence of even or odd approximants. For
(\ref{R3}), the behavior is even more interesting: it diverges for
$|q|<1$, but its sequence of approximants converge to different
values depending on the residue class modulo $3$ from which the
approximants are chosen. We show that (\ref{RR}) and (\ref{R3})
are part of the same phenomenon-- a phenomenon that we thoroughly
explore here.  We present a unified theory showing that this
behavior is typical of a larger class of continued fractions which
have multiple limits. We show that there is nothing special about
two or three limits, and that continued fractions with $m\ge 2$
limits arise just as naturally.

In the course of our investigations it became necessary to
investigate certain infinite products of matrices. These, in turn,
led to theorems about the limiting behavior of Poincar{\'e}-type
recurrences.  We obtain a theorem similar to that of  Perron on
the limiting behavior of Poincar{\'e}-type recurrence sequences in
the case where the eigenvalues are distinct roots of unity, but
our theorem gives more information. We begin by describing this
work.

\subsection{Poincar{\'e}-type recurrences}

Let the sequence $\{x_{n}\}_{n \geq 0}$ have the initial values
$x_{0}$, $\dots$, $x_{p-1}$ and be defined for larger values by
\begin{equation}\label{pereq}
x_{n+p}=\sum_{r=0}^{p-1}a_{n,r}x_{n+r},
\end{equation}
for $n \geq 0$. Suppose also that there are numbers $a_{0},\dots, a_{p-1}$ such that
\begin{align}
\label{asumineq}
&\lim_{n \to \infty}a_{n,r}=a_{r},& &0 \leq r \leq p-1.&
\end{align}

A recurrence of the form (\ref{pereq}) satisfying the condition (\ref{asumineq}) is
called a Poincar{\'e}-type recurrence, (\ref{asumineq}) being known as the
Poincar{\'e} condition. Such recurrences were initially studied by Poincar{\'e}
who proved that if the roots of the characteristic equation
\begin{equation}\label{chareq}
t^{p}-a_{p-1}t^{p-1}-a_{p-2}t^{p-2}- \dots -a_{0}=0
\end{equation}
have distinct norms, then the ratios of consecutive terms in the
recurrence (for any set of initial conditions) tend to one of the
roots.  See \cite{Poi}.  Because the roots are also the
eigenvalues of the associated companion matrix, they  are also
referred to as the eigenvalues of (\ref{pereq}).  This result was
improved by O. Perron, who obtained a number of theorems about the
limiting asymptotics of such recurrence sequences. Perron
\cite{Perron1} made a significant advance in 1921 when he proved
the following theorem which for the first time treated cases of
eigenvalues which repeat or are of equal norm.

\begin{theorem}\label{tp1}
Let the sequence $\{x_{n}\}_{n \geq 0}$ be defined by initial values
$x_{0}$, $\dots$, $x_{p-1}$ and by (\ref{pereq})
for $n \geq 0$. Suppose also that there are numbers $a_{0},\dots, a_{p-1}$ satisfying (\ref{asumineq}).
Let $q_{1},\,q_{2}, \dots q_{\sigma}$ be the distinct moduli of the characteristic equation (\ref{chareq})
and let $l_{\lambda}$ be the number of roots whose modulus is $q_{\lambda}$, multiple roots counted
according to multiplicity, so that
\[
l_{1}+l_{2}+\dots l_{\sigma}=p.
\]
Then, provided $a_{n,0}$ be different from zero for $n \geq 0$, the difference equation
\eqref{pereq} has a fundamental system of solutions,
which fall into $\sigma$ classes, such that, for the solutions of the $\lambda$-th class
and their linear combinations,
\[
\limsup_{n \to \infty}  \sqrt[n]{|x_{n}|}=q_{\lambda}.
\]
The number of solutions of the $\lambda$-th class is $l_{\lambda}$.
\end{theorem}

Thus when all of the characteristic roots have norm $1$, this theorem
gives that
\[
\limsup_{n  \to \infty}  \sqrt[n]{|x_{n}|}= 1.
\]

Another related paper is \cite{KL} where the authors study products
of matrices and give a sufficient condition for their boundedness. This
is then used to study ``equimodular'' limit periodic continued fractions, which
are limit periodic continued fractions in which the characteristic roots of the associated
$2 \times 2$ matrices are all equal in modulus. The matrix theorem in \cite{KL} can 
also be used to obtain results about the boundedness of recurrence sequences.
We study a more specialized situation
here and obtain far more detailed information as a consequence.

Our focus is on the case where the characteristic roots
are distinct roots of unity. Under a stronger condition 
than (\ref{asumineq}) we will show that all non-trivial solutions of such recurrences approach
a finite number of limits in a precisely controlled way. Specifically, our theorem
is:
\begin{theorem}\label{t2}
Let the sequence $\{x_{n}\}_{n \geq 0}$ be defined by initial values
$x_{0}$, $\dots$, $x_{p-1}$ and by (\ref{pereq})
for $n \geq 0$. Suppose also that there are numbers $a_{0},\dots, a_{p-1}$ such that
\begin{align*}
%\label{asumineq}
&\sum_{n=0}^{\infty}|a_{r}-a_{n,r}|<\infty,& &0 \leq r \leq p-1.&
\end{align*}
Suppose the solutions of (\ref{chareq})
are distinct roots of unity, say $\alpha_{1}$, $\dots$, $\alpha_{p}$.
Let $m$ be the least positive integer such that, for all $j \in \{0,1,\dots,p-1\}$,
$\alpha_{j}^{m} = 1$.
 Then, for $0 \leq j \leq m-1$, the subsequence
$\{x_{mn+j}\}_{n=0}^{\infty}$ converges. Set $l_j=\lim_{n\to\infty} x_{nm+j}$, for integers $j\ge 0$.
Then the (periodic) sequence $\{l_j\}$ satisfies the recurrence relation
\begin{equation}\label{limrecur}
l_{n+p}=\sum_{r=0}^{p-1}a_{r}l_{n+r},
\end{equation}
and thus there exist constants $c_1,\cdots,c_p$ such that
\begin{equation}\label{limform}
l_n=\sum_{i=1}^p c_i\alpha_i^n.
\end{equation}
\end{theorem}

Remark. In \cite{ABSYZ02} the authors study a recurrence which is of Poincar{\'e} type
and has $6$ limits. In section 4 we obtain this result of a special case of one
of our corollaries.   

Theorem \ref{t2} follows easily from Proposition \ref{tm} in section 2 which proves the
convergence of infinite products of matrices when the limits are taken in
arithmetic progressions. Proposition \ref{tm} is also our key to giving a unifying
theory of certain classes of continued fractions with multiple limits.

\subsection{Continued fractions with multiple limits}

Our general theorem on continued fractions with multiple limits is the following which
includes the multiple limit convergence behavior of both (\ref{RR}) and (\ref{R3}) as special cases.
\begin{theorem}\label{T1}
Let $\{p_{n}\}_{n\geq 1}$, $\{q_{n}\}_{n \ge 1}$ be
complex sequences satisfying
\begin{align*}
 &\sum_{n=1}^{\infty}|p_{n}|<\infty,& &\sum_{n=1}^{\infty}|q_{n}|<\infty.&
\end{align*}
Let $\omega_{1}$  and $\omega_{2}$ be distinct roots of unity
and let $m$ be the least positive integer such that $\omega_{1}^m=\omega_{2}^{m}=1$ .   Define
{\allowdisplaybreaks
\begin{equation*}
G:=
\frac{- \omega_{1}\omega_{2}+q_{1}}{\omega_{1}+\omega_{2}+p_{1}}
\+
\frac{-\omega_{1}\omega_{2}+q_{2}}{\omega_{1}+\omega_{2}+p_{2}}
\+
\frac{-\omega_{1}\omega_{2}+q_{3}}{\omega_{1}+\omega_{2}+p_{3}}
\+
\cds.
\end{equation*}
}
Let $\{P_{n}/Q_{n}/\}_{n=1}^{\infty}$ denote the sequence of approximants of $G$.
If $q_{n}\not = \omega_{1}\omega_{2}$ for any $n \geq 1$,
then $G$ does not converge. However, the sequences of numerators and denominators
in each of the $m$ arithmetic progressions modulo $m$ do converge. More precisely,
 there exist complex numbers $A_{0}, \dots  , A_{m-1}$ and $B_{0}, \dots  , B_{m-1}$
such that, for $0 \leq i< m$,
{\allowdisplaybreaks
\begin{align}\label{ABlim}
&\lim_{k \to \infty} P_{m\,k+i}=A_{i}, & &\lim_{k \to \infty} Q_{m\,k+i}=B_{i}.&
\end{align}
Extend the sequences $\{A_i\}$ and $\{B_i\}$ over all integers by making them
periodic modulo $m$ so that (\ref{ABlim}) continues to hold. Then for integers $i$,
\begin{equation}\label{lol}
A_i=\left(\frac{A_1-\omega_2A_0}{\omega_1-\omega_2}\right)\omega_1^i
+\left(\frac{\omega_1A_0-A_1}{\omega_1-\omega_2}\right)\omega_2^i,
\end{equation}
and 
\begin{equation}\label{lolb}
B_i=\left(\frac{B_1-\omega_2B_0}{\omega_1-\omega_2}\right)\omega_1^i
+\left(\frac{\omega_1B_0-B_1}{\omega_1-\omega_2}\right)\omega_2^i.
\end{equation}

Moreover, 
\begin{equation}\label{detAB}
A_iB_{j}-A_{j}B_i=-(\omega_1\omega_2)^{j+1}
\frac{\omega_1^{i-j}-\omega_2^{i-j}}{\omega_1-\omega_2}
\prod_{n=1}^{\infty}\left(1-\frac{q_{n}}{\omega_1\omega_2}\right).
\end{equation}
}
Put  $\omega_{1} := \exp( 2 \pi i a/m)$,  $\omega_{2} := \exp (2 \pi i b/m)$, $0 \leq a <b <m$, 
and $r:=m/\gcd(b-a,m)$. Then $G$ has $r$ distinct limits  
in $\hat{\mathbb{C}}$ which are given by $A_j/B_j$, $1\le j\le r$. Finally, for $k \geq 0$ and $1 \leq j \leq r$,  
{\allowdisplaybreaks
\begin{equation*}
\frac{A_{j+kr}}{B_{j+kr}}=\frac{A_{j}}{B_{j}}.
\end{equation*}

}
\end{theorem}

To see how the behavior of (\ref{RR}) for $|q|>1$ follows from this, observe that
by the standard equivalence transformation for continued fractions, (\ref{RR})
is equal to
 \[
 1+
\frac{1}{1/q}
\+
\frac{1}{1/q}
\+
\frac{1}{1/q^{2}}
 \+
\frac{1}{1/q^{2}}
\cds
\+
\frac{1}{1/q^{n}}
 \+
\frac{1}{1/q^{n}}
\cds
.
\]
We can now apply Theorem \ref{T1} with $\omega_{1}=-1$,  $\omega_{2}=1$
(so $m=2$), $q_{n}=0$ and $p_{2n-1}=p_{2n}=1/q^{n}$ to get that the continued fraction
does not converge, but that the sequence of approximants in each of the arithmetic
progressions modulo $2$ do converge.

The behavior of (\ref{R3}) is similarly a special case. Just let $\omega_{1}=\exp(2 \pi i/6)$,  $\omega_{2}=\exp(-2 \pi i/6)$
$($so that $\omega_{1}+\omega_{2}=\omega_{1}\omega_{2}=1$$)$, $g_n=0$ and $f_n=q^n$.
Theorem \ref{T1} then gives that (\ref{R3}) has three limits for $|q|<1$. 

We refer to the number $r$ in the theorem as the {\it rank} of 
the
continued fraction. 
A remarkable consequence of this theorem is that because of (\ref{lol}) and (\ref{lolb}),
to compute all the limits of a continued
fraction with of rank $r$, one only needs to know the first two limits of the numerator
and denominator convergents. In fact, taking $i=0$ and $j=1$ in (\ref{detAB}), it can be seen 
that one only needs to know the value of three of the four limits $\{A_0, A_1, B_0, B_1\}$.

Another interesting consequence of Theorem \ref{T1} is that the fundamental Stern-Stolz
divergence theorem \cite{LW92}  is an immediate corollary. In fact, the Stern-Stolz theorem will be found
to be the beginning of an infinite sequence of similar theorems all of which are special cases
of our theorem. See Corollaries  \ref{Stern-Stolz}-\ref{cor1} and Example \ref{SSseq}.
These  consequences of Theorem \ref{T1} are explored after its proof.

\subsection{A generalization of the Ramanujan continued fraction with three limits}

In a recent paper \cite{ABSYZ02}, the authors proved a claim
made
by Ramanujan in his lost notebook (\cite{S88}, p.45) about (\ref{R3}). To
describe Ramanujan's claim, we first need some notation.
 Throughout take $q\in\mathbb{C}$
with $|q|<1$. The following standard notation for $q$-products
will also be employed:
\begin{align*}
&(a)_{0}:=(a;q)_{0}:=1,& &
(a)_{n}:=(a;q)_{n}:=\prod_{k=0}^{n-1}(1-a\,q^{k}),& & \text{ if }
n \geq 1,&
\end{align*}
and
\begin{align*}
&(a;q)_{\infty}:=\prod_{k=0}^{\infty}(1-a\,q^{k}),& & |q|<1.&
\end{align*}
Set $\omega = e^{2 \pi i/3}$. Ramanujan's claim was that, for
$|q|<1$,
{\allowdisplaybreaks
\small{
\begin{equation}\label{3lim1}
\lim_{n \to \infty}
 \left (
\frac{1}{1}
 \-
\frac{1}{1+q}
\-
\frac{1}{1+q^2}
\-
\cds
\-
\frac{1}{1+q^n +a}
\right )
=
-\omega^{2}
\left (
\frac{\Omega - \omega^{n+1}}{\Omega - \omega^{n-1}}
\right ).
\frac{(q^{2};q^{3})_{\infty}}{(q;q^{3})_{\infty}},
\end{equation}
}
}
where
{\allowdisplaybreaks
\begin{equation*}
\Omega :=\frac{1-a\omega^{2}}{1-a \omega}
\frac{(\omega^{2}q,q)_{\infty}}{(\omega q,q)_{\infty}}.
\end{equation*}
} Ramanujan's notation is confusing, but what his claim means is
that the limit exists as $n \to \infty$ in each of the three
congruence classes modulo 3, and that the limit is given by the
expression on the right side of (\ref{3lim1}). Also, the appearance
of the variable $a$ in this formula is a bit of a red herring; from elementary
properties of continued fractions, one can derive the result for
general $a$ from the $a=0$ case. 

Here we examine a direct generalization of Ramanujan's continued
fraction which has $k$ limits, for an arbitrary positive integer
$k \geq 2$, and evaluate these limits in terms of ratios of
certain unusual $q$-series. Let $m$ be any arbitrary
integer greater than $2$, let $\omega$ be a primitive $m$-th root
of unity and, for ease of notation,
  let $\bar{ \omega} = 1/\omega$.  Define
{\allowdisplaybreaks
\begin{equation*}
G(q):= \frac{1}{1} \- \frac{1}{\omega + \bar{ \omega}+q} \-
\frac{1}{\omega + \bar{ \omega}+q^2} \- \frac{1}{\omega + \bar{
\omega}+q^3} \+ \cds.
\end{equation*}
} For  $j \geq 0$, $k \geq 0$  and $ i \in \mathbb{Z}$ define
%{\allowdisplaybreaks
\begin{equation*}
F_{k}(\omega, i,j,q) :=  \sum_{u=0}^{ \lfloor \frac{mk+i}{2}
\rfloor '} (q^{u+1};q)_{j}(\omega^{2u-i} + \omega^{i-2u}),
\end{equation*}
%}
where the summation $\sum _{u=0}^{ \lfloor \frac{mk+i}{2} \rfloor
'}$ means that if $mk+i$ is even, the final term in the sum
%, corresponding to
%$u=\lfloor (mk+i)/2 \rfloor $,
is $(q^{(mk+i)/2+1};q)_{j}$, rather
 than
$2(q^{(mk+i)/2+1};q)_{j}$. Then $F(\omega, i,j,q):=\lim_{k \to
\infty}F_{k}(\omega, i,j,q)$ exists and is finite. We prove the
following theorem.
\begin{theorem}\label{ramgentheor}
Let $\omega$ be a primitive $m$-th root of unity and
  let $\bar{ \omega} = 1/\omega$.
Let $1\leq i \leq m$. Then
{\allowdisplaybreaks
\begin{multline}\label{ramgentheoreq}
\lim_{k \to \infty}
\frac{1}{\omega + \bar{ \omega}+q}
\-
\frac{1}{\omega + \bar{ \omega}+q^2}
\-
\cds
\frac{1}{\omega + \bar{ \omega}+q^{mk+i}}\\
= \frac{\sum_{j=0}^{\infty}\frac{q^{j(j+3)/2}}{(q;\,q)_{j}^{2}}F(\omega,i-j-2,j,q)}
{\sum_{j=0}^{\infty}\frac{q^{j(j+1)/2}}{(q;\,q)_{j}^{2}}F(\omega,i-j-1,j,q)}.
\end{multline}
}
Moreover, the continued fraction has rank $m$ when $m$ is odd, and 
rank $m/2$ when $m$ is even.
\end{theorem}
We state the result for the first tail of $G(q)$, rather than
$G(q)$ itself, for aesthetic reasons.

\subsection{Analytic continued fractions with multiple limits, via Daniel Bernoulli's Theorem}

In  \cite{ABSYZ02}, the authors also describe a general class of
analytic continued fractions with three limit points.  Our final
results are to give an alternative derivation, using Daniel
Bernoulli's continued fraction, of this class of analytic
continued fractions, and to generalize this class, again using
Bernoulli's continued fraction, by showing how to construct
analytic continued fractions with $m$ limit points, for an
arbitrary positive integer $m \geq 2$. Recall that Bernoulli's
result is an easy formula for a continued fraction whose sequence
of approximants agrees exactly with any prescribed sequence (See
Proposition \ref{pber}). If the original sequence is of a simple
kind, then Bernoulli's continued fraction offers no advantage over
the original sequence and so is in a certain sense ``trivial''.
(It is just an obscure way of writing down a simple sequence.) Our
results in Theorems 1 and 2 are deeper and do not arise in this
way. Using Bernoulli's continued fraction we prove the following
theorem which gives a general class of continued fractions having
$m$ limit points. The sequence equal to the $n$'th approximant of
this continued fraction is constructed from taking the union of
$m$ sequences. This gives rise to the $m$ limits. We have included
this result to put into perspective the special case of
Bernoulli's formula that is given in \cite{ABSYZ02}.

\begin{theorem}\label{t4}
Let $G(z)$ be analytic in the closed unit disc and suppose
{\allowdisplaybreaks
\begin{align}\label{gcon}
&|G(z)|< 1/2,&  &   |z|<1.&
\end{align}
}
Let $m$ be a positive integer, $m \geq 2$.
Define
{\allowdisplaybreaks
\begin{equation} \gamma_{n} =
\begin{cases} 1-m, \,\,\,\,n \equiv 0 (\text{mod } m ),\\
1, \text{\phantom{sdasda}otherwise}.
\end{cases}
\end{equation}
}
 Then, for $|z|<1$, the continued fraction
{\allowdisplaybreaks
\begin{multline}\label{bernoucf}
\frac{G(z^{2})-G(z)+1}{1}
\-
\frac{G(z^{3})-G(z^{2})+1}{G(z^{3})-G(z)+2}\\
\+
K_{n=3}^{\infty}\frac{-(G(z^{n+1})-G(z^{n})+\gamma_{n})(G(z^{n-1})-G(z^{n-2})+\gamma_{n-2})}
{G(z^{n+1})-G(z^{n-1})+\gamma_{n}+\gamma_{n-1}}.
\end{multline}
}
 has exactly  the $m$ limits  $G(0)-G(z)+i$, $0 \leq i \leq m-1$.
\end{theorem}

\section{A Result on Infinite Matrix Products}
The convergence results of this paper follow from the following proposition.
\begin{proposition}\label{tm}
Let $p\geq 2$ be an integer and let $M$ be a $p \times p$ matrix that is diagonalizable
and whose eigenvalues are roots of unity. Let $I$ denote the $p \times p$ identity matrix and let
$m$ be the least positive integer such that
\[
M^{m}=I.
\]
For a $p \times p$ matrix $G$, let
\[
||G||_{\infty}= \max_{1\leq i,j\leq p}|G^{(i,j)}|.
\]
Suppose $\{D_{n}\}_{n=1}^{\infty}$ is a sequence of matrices such that
\[
\sum_{n=1}^{\infty}||D_{n}-M||_{\infty}< \infty.
\]
 Then
\[F:=\lim_{k \to \infty} \prod_{n=1}^{km}D_{n}
\]
 exists. Here the matrix product means either $D_{1}D_{2} \dots$ or $\dots D_{2}D_{1}$.
Further, for each $j$, $0 \leq j \leq m-1$,
\[
\lim_{k \to \infty} \prod_{n=1}^{km+j}D_{n} = M^{j}F \text{ or } FM^{j},
\]
depending on whether the products are taken to the left or right.
\end{proposition}
We prove the proposition for the products $D_{1}D_{2} \dots$ only,
since the other case is virtually identical.
We need two preliminary lemmas.
\begin{lemma}
For $n \geq 0$, define
{\allowdisplaybreaks
\[
U_{n}=\prod_{j=1}^{m}D_{mn+j}.
\]
}
Then there exists  a sequence $\{\epsilon_{n}\}$ with
$\sum_{n=0}^{\infty}\epsilon_{n}<\infty$ and an absolute constant $A$ such that
\[
||U_{n}-I||_{\infty} \leq \epsilon_{n}A.
\]
\end{lemma}
\begin{proof}
Let $\epsilon_{n}=\max_{1\leq j \leq m}||D_{mn+j}-M||_{\infty}$. Define $E_{mn+j}$ by
\[
D_{mn+j}=M+\epsilon_{n}E_{mn+j}.
\]
(If $\epsilon_{n}=0$, define $E_{mn+j}$ to be the $p\times p$ zero matrix). Note that the entries in
each matrix  $E_{mn+j}$ are bounded in absolute value by 1. Let the matrix $R_{n}$ be as defined below.
\begin{align*}
U_{n}= \prod_{j=1}^{m}D_{mn+j}&=\prod_{j=1}^{m}\left ( M + \epsilon_{n}E_{mn+j} \right)\\
&:=M^{m}+\epsilon_{n}R_{n}=I+\epsilon_{n} R_{n}.
\end{align*}
The elements of all the matrices $R_{n}$ for $n \geq 0$ are absolutely bounded
(independent of $n$) since $R_{n}$ is formed from
a sum of at most $2^{m}$ products of matrices, where each product contains $m$ matrices
and the entries in each matrix are bounded by $\max \{||M||_{\infty},\epsilon_{n} \}$. Let
$A= \sup \{||R_{n}||_{\infty}\}$. Then
\begin{equation*}
||U_{n}-I||_{\infty} = \epsilon_{n}||R_{n}||_{\infty} \leq \epsilon_{n}A.
\end{equation*}
\end{proof}

\begin{lemma}
With the notation of the previous lemma, define
\[
F_{r}=\prod_{n=0}^{r}U_{n}.
\]
Then $\lim_{r \to \infty}F_{r}$ exists.
\end{lemma}

\begin{proof}
 Let $A$ be as defined in the previous lemma. \\
\emph{Claim 1:}
\[
||F_{r}||_{\infty}\leq \prod_{j=0}^{r}\left (1+p\epsilon_{j}A\right).
\]
\emph{Proof of claim.}. For $r=0$, $F_{0}=U_{0}$ and
\[
|F_{0}^{(i,j)}|=|U_{0}^{(i,j)}|\leq 1+\epsilon_{0}A \leq 1+p\epsilon_{0}A .
\]
Assume the claim is true for $r=0,1,\dots s$.
\begin{align*}
|F_{s+1}^{(i,j)}|&= \left | \sum_{k=1}^{p} F_{s}^{(i,k)}U_{s+1}^{(k,j)} \right | \leq
 \sum_{k=1}^{p} \left |F_{s}^{(i,k)}\right | \left |U_{s+1}^{(k,j)} \right |
\leq \prod_{j=0}^{s}(1+p\epsilon_{j}A)  \sum_{k=1}^{p}\left |U_{s+1}^{(k,j)} \right | \\
&\leq \prod_{j=0}^{s}(1+p\epsilon_{j}A)  (1+p \epsilon_{s+1}A)
= \prod_{j=0}^{s+1}(1+p\epsilon_{j}A).
\end{align*}
In particular, note that, for each $r \geq 0$ and each index $(i,j)$,
\[
|F_{r}^{(i,j)} | \leq \prod_{j=0}^{\infty}(1+p\epsilon_{j}A):=C.
\]
Note that the infinite product converges, since
$\sum_{n=0}^{\infty}\epsilon_{n} \leq \sum_{n=0}^{\infty}||D_{n}-M|| < \infty$.

\emph{Claim 2:} For each index $(i,j)$, the sequence
$\{F_{r}^{(i,j)}\}$ is Cauchy, and hence convergent.\\
\emph{Proof of claim.}  By definition,
\[
F_{r+1}-F_{r} = F_{r}U_{r+1}-F_{r} = F_{r}(U_{r+1}-I).
\]
Hence,
\begin{align*}
|F_{r+1}^{(i,j)}-F_{r}^{(i,j)}| &=
\left | \sum_{k=1}^{p}F_{r}^{(i,k)}(U_{r+1}-I)^{(k,j)} \right |\\
&\leq \sum_{k=1}^{p}|F_{r}^{(i,k)}|\left |(U_{r+1}-I)^{(k,j)} \right |
\leq pCA\epsilon_{r+1},
\end{align*}
where $C$ is as defined above.
This is sufficient to show the sequence is Cauchy, since
$\sum_{n=0}^{\infty}\epsilon_{n}<0$.
Define the matrix $F$ by
\[
F^{(i,j)}:=\lim_{r \to \infty}F_{r}^{(i,j)}.
\]
\end{proof}
\emph{Proof of Proposition \ref{tm}.} This now follows easily from the above lemma, since
\begin{align*}
\lim_{k \to \infty} \prod_{n=1}^{km+j}D_{n} &=
\lim_{k \to \infty} \prod_{n=1}^{km}D_{n}\prod_{n=km+1}^{km+j}D_{n}\\
&=\lim_{k \to \infty} F_{k-1}\prod_{n=km+1}^{km+j}D_{n}
=FM^{j}.
\end{align*}

\begin{flushright}
$ \Box $
\end{flushright}

\section{Continued Fractions With Multiple Limits}
Proposition \ref{tm} allows us to construct non-trivial divergent continued fractions
whose sequences of approximants in each of the arithmetic progressions
modulo $m$ converge. We now prove Theorem \ref{T1}.

\begin{proof}[Proof of Theorem \ref{T1}]
Let
\[
M=\left (
\begin{matrix}
&\omega_{1}+\omega_{2} &1 \\
&-\omega_{1}\omega_{2} &0
\end{matrix}
\right ).
\]
It follows easily from the identity
{\allowdisplaybreaks
\begin{align*}
\left (
\begin{matrix}
&1 &
1\\
&\phantom{as} & \phantom{as} \\
&- \omega_{2}& -\omega_{1}
\end{matrix}
\right )
\left (
\begin{matrix}
&
\omega_{1}
&0 \\
&\phantom{as} & \phantom{as} \\
&0 &\omega_{2}
\end{matrix}
\right )
\left (
\begin{matrix}
&1 &
1\\
&\phantom{as} & \phantom{as} \\
&- \omega_{2}& -\omega_{1}
\end{matrix}
\right )^{-1}
=
\left (
\begin{matrix}
&\omega_{1}+\omega_{2} & 1 \\
&\phantom{as} & \phantom{as} \\
&-\omega_{1}\omega_{2} & 0
\end{matrix}
\right )
\end{align*}
}
that
\begin{equation}\label{mpower}
M^{j}=
\left (
\begin{matrix}
\displaystyle{
 \frac{{{\omega_1}}^{1 + j} - {{\omega_2}}^{1 + j}}
   {{\omega_1} - {\omega_2}} }
& \displaystyle{
\frac{{{\omega_1}}^j - {{\omega_2}}^j}
   {{\omega_1} - {\omega_2}}} \\ 
\phantom{as} & \phantom{as} \\
- \displaystyle{
\frac{{\omega_1}\,{\omega_2}\,
\left( {{\omega_1}}^j - {{\omega_2}}^j \right) }{{\omega_1} - {\omega_2}}}
       
& \displaystyle{
\frac{- {{\omega_1}}^j\,{\omega_2}   + 
     {\omega_1}\,{{\omega_2}}^j}{{\omega_1} - {\omega_2}} }
\end{matrix}
\right ),
\end{equation}
and thus that
{\allowdisplaybreaks
\begin{align*}
&M^{m}=\left (
\begin{matrix}
&1 &0 \\
&0 &1
\end{matrix}
\right ),&
&M^{j} \not = \left (
\begin{matrix}
&1 &0 \\
&0 &1
\end{matrix}
\right ),& &1 \leq j < m.&
\end{align*}
}
 For $n \geq 1$, define
{\allowdisplaybreaks
\begin{equation*}
D_{n} :=
\left (
\begin{matrix}
&\omega_{1}+\omega_{2} +p_{n} & 1 \\
&\phantom{as} & \phantom{as} \\
&-\omega_{1}\omega_{2} +q_{n}& 0
\end{matrix}
\right ).
\end{equation*}
}
Then
\[
\sum_{n \ge 1}||D_{n}-M||_{\infty}<\infty.
\]
Further,
\[
||D_{n}-M||_{\infty} = \max \{ |p_{n}|, |q_{n}|\}.
\]
 Thus the matrix $M$ and the matrices
$D_{n}$ satisfy the conditions of Proposition \ref{tm}.
 Let the matrices $F_{i}$ and $F$ have the same meaning as in the proof of Proposition \ref{tm}.

By the correspondence between matrices and continued fractions,
{\allowdisplaybreaks
\begin{align}\label{ABmat}
\left (
\begin{matrix}
& P_{mn+i} &P_{mn+i-1}  \\
&\phantom{as} & \phantom{as} \\
& Q_{mn+i} & Q_{mn+i-1}
\end{matrix}
\right )=
\left (
\begin{matrix}
& 0 & 1 \\
&\phantom{as} & \phantom{as} \\
& 1& 0\end{matrix}
\right )
\prod_{j=1}^{mn+i}
D_{j} &\\
\phantom{as}& \notag \\
=
\left (
\begin{matrix}
& 0 & 1 \\
&\phantom{as} & \phantom{as} \\
& 1& 0\end{matrix}
\right )
F_{n-1}\prod_{j=mn+1}^{mn+i}
D_{j}.& \notag
\end{align}
} Now let $n \to \infty$ to get that
{\allowdisplaybreaks
\begin{equation}\label{abproof}
\lim_{n \to \infty}
\left (
\begin{matrix}
& P_{mn+i} &P_{mn+i-1} \\
&\phantom{as} & \phantom{as} \\
& Q_{mn+i} & Q_{mn+i-1}
\end{matrix}
\right )
=
\left (
\begin{matrix}
& 0 & 1 \\
&\phantom{as} & \phantom{as} \\
& 1& 0\end{matrix}
\right )
F\,M^{i}.
\end{equation}
}
This proves \eqref{ABlim}. 

Now let $A_i:=\lim_{n\to\infty}P_{mn+i}$, and 
$B_i:=\lim_{n\to\infty}Q_{mn+i}$. Notice by definition that the sequences $\{A_i\}$
and $\{B_i\}$ are periodic modulo $m$.

It easily follows from  (\ref{abproof}) that
{\allowdisplaybreaks
\begin{equation*}
\left (
\begin{matrix}
& A_{i} &A_{i-1} \\
&\phantom{as} & \phantom{as} \\
& B_{i} & B_{i-1}
\end{matrix}
\right )
=
\left (
\begin{matrix}
& A_{j} &A_{j-1} \\
&\phantom{as} & \phantom{as} \\
& B_{j} & B_{j-1}
\end{matrix}
\right )
M^{i-j}.
\end{equation*}
}
(\ref{mpower}) also gives that
{\allowdisplaybreaks
\begin{equation}\label{2da}
A_{i} = A_{j}
\displaystyle{
 \frac{{{\omega_1}}^{1 + i-j} - {{\omega_2}}^{1 +i- j}}
   {{\omega_1} - {\omega_2}} }
-A_{j-1}
\displaystyle{
\frac{{\omega_1}\,{\omega_2}\,
       \left( {{\omega_1}}^{i-j} - {{\omega_2}}^{i-j} \right) }{{\omega_1} - {\omega_2}}},
\end{equation}
and
\begin{equation}\label{2db}
B_{i} = B_{j}
\displaystyle{
 \frac{{{\omega_1}}^{1 + i-j} - {{\omega_2}}^{1 +i- j}}
   {{\omega_1} - {\omega_2}} }
-B_{j-1}
\displaystyle{
\frac{{\omega_1}\,{\omega_2}\,
       \left( {{\omega_1}}^{i-j} - {{\omega_2}}^{i-j} \right) }{{\omega_1} - {\omega_2}}}.
\end{equation}
}
Thus
\begin{equation*} 
A_{i}B_{j}-A_{j}B_{i}=
 \frac{\left( {A_j}\,{B_{-1 + j}} - 
        {A_{-1 + j}}\,{B_j} \right) \,{\omega_1}\,{\omega_2}\,
      \left( {{\omega_1}}^{i - j} - {{\omega_2}}^{i - j} \right) }
      {{\omega_1} - {\omega_2}}.
\end{equation*}
(\ref{lol}) and (\ref{lolb}) follow from (\ref{2da}) and (\ref{2db}) by setting $j=1$.
(\ref{detAB}) follows after applying the determinant formula
\begin{align*}
A_{j}B_{j-1}-A_{j-1}B_{j} 
&= -\lim_{k \to \infty} \prod_{n=1}^{mk+j}(\omega_{1}\omega_{2}-q_{n})\\
&= - (\omega_{1}\omega_{2})^{j}\prod_{n=1}^{\infty}
\left ( 1-\frac{q_{n}}{\omega_{1}\omega_{2}} \right ).
\end{align*}
Since $\sum_{j=1}^{\infty}|q_{j}|$ converges to a finite value,
 the infinite product on the right side converges.  

For the continued fraction to converge, $A_{i}B_{i-1}-A_{i-1}B_{i}=0$ is required. However,
(\ref{detAB}) shows that this is not the case.

Also from \eqref{detAB} we have that $A_{i}B_{j}-A_{j}B_{i}=0$, and thus that 
$A_{i}/B_{i}=A_{j}/B_{j}$ in the extended complex plane, if and only if
\[
\omega_{1}^{i-j}=\omega_{2}^{i-j}.
\]
This happens if and only if 
\[
(i-j)a \equiv (i-j) b \mod{m}.
\]
It follows easily from this that 
\begin{align*}
&\frac{A_{j}}{B_{j}},& &  1 \leq j \leq \frac{m}{\gcd(b-a,m)}=:r,&
\end{align*}
 are distinct and 
that 
{\allowdisplaybreaks
\begin{align*}
&\frac{A_{j+kr}}{B_{j+kr}}=\frac{A_{j}}{B_{j}},&& 1 \leq j \leq r, &&k \geq 0.&
\end{align*}
}
\end{proof}

It is easy to derive  general divergence results from this theorem, including the classical
Stern-Stolz Theorem \cite{LW92}. In fact, Stern-Stolz can be seen as the beginning of an infinite family
of divergence theorems. We first derive the Stern-Stolz theorem as a corollary, generalize 
it, then give a corollary describing the infinite family. Last, we  list the first
few examples in the infinite family.

\begin{corollary} (Stern-Stolz)\label{Stern-Stolz}
Let the sequence $\{b_{n}\}$ satisfy $\sum |b_{n}| <\infty$. Then the continued fraction
\[
b_{0}+K_{n=1}^{\infty}\frac{1}{b_{n}}
\]
diverges. In fact, for $p=0,1$,
\begin{align*}
&\lim_{n \to \infty}P_{2n+p}=A_{p} \not = \infty,& &\lim_{n \to \infty}Q_{2n+p}=B_{p} \not = \infty,&
\end{align*}
and
\[
A_{1}B_{0}-A_{0}B_{1} = 1.
\]
\end{corollary}
\begin{proof}
This follows immediately from Theorem \ref{T1}, upon setting $\omega_{1}=1$,
$\omega_{2} = -1$ (so $m=2$), $q_{n}=0$ and $p_{n} = b_{n}$.
\end{proof}

Note that by taking $p_n=a_n$ we immediately obtain a generalization.

\begin{corollary}\label{SSgen}
Let the sequences $\{a_n\}$ and $\{b_{n}\}$ satisfy $a_n\ne -1$ for $n\ge 1$,  $\sum |a_{n}| <\infty$ and
$\sum |b_{n}| <\infty$.
Then the continued fraction
\[
b_{0}+K_{n=1}^{\infty}\frac{1+a_n}{b_{n}}
\]
diverges. In fact, for $p=0,1$,
\begin{align*}
&\lim_{n \to \infty}P_{2n+p}=A_{p} \not = \infty,& &\lim_{n \to \infty}Q_{2n+p}=B_{p} \not = \infty,&
\end{align*}
and
\[
A_{1}B_{0}-A_{0}B_{1} = \prod_{n=1}^{\infty}(1+a_n).
\]
\end{corollary}
\begin{proof}
This follows immediately from Theorem \ref{T1}, upon setting $\omega_{1}=1$,
$\omega_{2} = -1$ (so $m=2$), $q_{n}=a_n$ and $p_{n} = b_{n}$.
\end{proof}

We have not been able to find Corollary \ref{SSgen} in the literature, but it is surly known
in principle.

The infinite family of Stern-Stolz type theorems is described by the following corollary.

\begin{corollary}\label{cor1}
Let the sequences $\{a_n\}$ and $\{b_{n}\}$ satisfy
$a_{n} \not = 1$ for $n\ge 1$, $\sum |a_{n}| <\infty$ and $\sum |b_{n}| <\infty$. Let 
$m \geq 3$ and let
$\omega_{1}$ be a primitive $m$-th root of unity. Then the continued fraction 
\[
b_{0}+K_{n=1}^{\infty}\frac{-1+ a_{n}}{\omega_{1}+ \omega_{1}^{-1}+b_{n}}
\]
does not converge, but the numerators and denominators in each of the
$m$ arithmetic progressions modulo $m$ do converge.
If $m$ is even, then for $i \leq p \leq m/2$,
\begin{align*}
&\lim_{n \to \infty}P_{mn+p}=-\lim_{n \to \infty}P_{mn+p+m/2}=A_{p} \not = \infty,& \\
&\lim_{n \to \infty}Q_{mn+p}=-\lim_{n \to \infty}Q_{mn+p+m/2}=B_{p} \not = \infty.&
\end{align*}
If $m$ is odd, then 
the continued fraction has rank $m$. If $m$ is even, then
the continued fraction has rank $m/2$. Further, for $2 \leq p \leq m'$, where 
$m'=m$ if $m$ is odd and $m/2$ if $m$ is even,
\[
A_{p}B_{p-1}-A_{p-1}B_{p} = -\prod_{n=1}^{\infty}(1-a_n).
\] 
\end{corollary}
\begin{proof}
In Theorem \ref{T1}, let $\omega_{2} = 1/\omega_{1}$.
\end{proof}
Some explicit examples are given below.

\begin{example}\label{SSseq}
Let the sequences $\{a_n\}$ and $\{b_{n}\}$ satisfy
$a_{n} \not = 1$ for $n\ge 1$, $\sum |a_{n}| <\infty$ and $\sum |b_{n}| <\infty$. Then
each of the following continued fractions diverges:

(i) The following continued fraction has rank three:
 \begin{equation}\label{c3} 
b_{0}+K_{n=1}^{\infty}\frac{-1+ a_{n}}{1+b_{n}}.
\end{equation}
 In fact, for $p=1,2,3$,
\begin{align*}
&\lim_{n \to \infty}P_{6n+p}=-\lim_{n \to \infty}P_{6n+p+3}=A_{p} \not = \infty,& \\
&\lim_{n \to \infty}Q_{6n+p}=-\lim_{n \to \infty}Q_{6n+p+3}=B_{p} \not = \infty.&
\end{align*}
(ii) The following continued fraction has rank four:
 \begin{equation}\label{c4} 
b_{0}+K_{n=1}^{\infty}\frac{-1+ a_{n}}{\sqrt{2}+b_{n}}.
\end{equation}
In fact, for $p=1,2,3,4$,
\begin{align*}
&\lim_{n \to \infty}P_{8n+p}=-\lim_{n \to \infty}P_{8n+p+4}=A_{p} \not = \infty,& \\
&\lim_{n \to \infty}Q_{8n+p}=-\lim_{n \to \infty}Q_{8n+p+4}=B_{p} \not = \infty.&
\end{align*}
(iii) The following continued fraction has rank five:
\begin{equation}\label{c5} 
b_{0}+K_{n=1}^{\infty}\frac{-1+ a_{n}}{(1-\sqrt{5})/2+b_{n}}.
\end{equation}
 In fact, for $p=1,2,3,4,5$,
\begin{align*}
&\lim_{n \to \infty}P_{5n+p}=A_{p} \not = \infty,& &\lim_{n \to \infty}Q_{5n+p}=B_{p} \not = \infty.&
\end{align*}
(iv) The following continued fraction has rank six:
 \begin{equation}\label{c6} 
b_{0}+K_{n=1}^{\infty}\frac{-1+ a_{n}}{\sqrt{3}+b_{n}}.
\end{equation}
In fact, for $p=1,2,3,4,5,6$,
\begin{align*}
&\lim_{n \to \infty}P_{12n+p}=-\lim_{n \to \infty}P_{12n+p+6}=A_{p} \not = \infty,& \\
&\lim_{n \to \infty}Q_{12n+p}=-\lim_{n \to \infty}Q_{12n+p+6}=B_{p} \not = \infty.&
\end{align*}
In each case we have, for $p$ in the appropriate range, that
\[
A_{p}B_{p-1}-A_{p-1}B_{p} = -\prod_{n=1}^{\infty}(1-a_n).
\]
\end{example}
\begin{proof}
In Corollary \ref{cor1}, set 

(i) $\omega_{1}=\exp ( 2 \pi i/6)$;

(ii) $\omega_{1} = \exp (2 \pi i/8)$;

(iii) $\omega_{1}=\exp ( 2 \pi i/5)$;

(iv)  $\omega_{1}=\exp ( 2 \pi i/12)$.
\end{proof}

The cases $\omega_{1}=\exp ( 2 \pi i/m)$, $m=3,4,10$ give a continued fractions 
that are the same as those above after an equivalence
transformation and renormalization of the sequences $\{a_n\}$ and$\{b_n\}$.
Note that the continued fractions (\ref{c4}) and (\ref{c6}) are, after an equivalence transformation
and renormalizing the sequences $\{a_n\}$ and $\{b_n\}$, of the forms
\[
b_{0}+K_{n=1}^{\infty}\frac{-2+ a_{n}}{{2}+b_{n}},
\]
and
\[
b_{0}+K_{n=1}^{\infty}\frac{-3+ a_{n}}{{3}+b_{n}},
\]
respectively.  Also, it should be mentioned that Theorem 3.3 of \cite{ABSYZ02} is essentially the
special case $a_n=0$ of part (i) of our example.

Theorem \ref{T1} now makes it trivial to construct $q$-continued fractions with arbitrarily many limits.

\begin{example}
Let $f(x)$, $g(x) \in \mathbb{Z}[q][x]$ be
polynomials with  zero constant term.
Let $\omega_{1}$, $\omega_{2}$ be distinct roots of unity
and suppose  $m$ is the least positive integer such that $\omega_{1}^m=\omega_{2}^{m}=1$ .
Define
{\allowdisplaybreaks
\begin{equation*}
G(q):=
\frac{-\omega_{1}\omega_{2}+g(q)}{\omega_{1}+\omega_{2}+f(q)}
\+
\frac{-\omega_{1}\omega_{2}+g(q^2)}{\omega_{1}+\omega_{2}+f(q^2)}
\+
\frac{-\omega_{1}\omega_{2}+g(q^3)}{\omega_{1}+\omega_{2}+f(q^3)}
\+
\cds.
\end{equation*}
}
Let $|q|<1$.
If $g(q^n) \not = \omega_{1}\omega_{2}$ for any $n \geq 1$,
then $G(q)$ does not converge.
However,  the sequences of approximants of
$G(q)$ in each of the $m$ arithmetic progressions modulo $m$ converge
to values in $\hat{\mathbb{C}}$. The continued fraction has rank $m/\gcd(b-a,m)$, 
where $a$ and $b$ are as defined in Theorem \ref{T1}. 
\end{example}

From this example we can conclude that (\ref{RR}) and (\ref{R3}) are far from unique
examples and many other $q$-continued fractions with multiple
limits can be immediately written down. For example, to Ramanujanize a bit, one can
immediately see that the continued fractions
\begin{equation}
\K{n\ge1}{\infty}\frac{-1/2}{1+q^n}\qquad\text{ and }\qquad \K{n\ge1}{\infty}\frac{-1/2+q^n}{1+q^n}
\end{equation}
both have rank four, while the continued fractions
\begin{equation}
\K{n\ge1}{\infty}\frac{-1/3}{1+q^n}\qquad\text{ and }\qquad \K{n\ge1}{\infty}\frac{-1/3+q^n}{1+q^n}
\end{equation}
both have rank six. We do not dwell further on $q$-continued fractions here, but 
in section 5 we will study a direct generalization of  (\ref{R3}).

\section{Recurrence relations with characteristic equations whose roots are roots of unity}

Theorem \ref{t2} follows easily from Proposition \ref{tm}. We now prove Theorem \ref{t2}.

\begin{proof}[Proof of Theorem \ref{t2}]
Define
\begin{equation*}
M:=
\left (
\begin{matrix}
&a_{p-1} &a_{p-2}& \dots & a_{1}&a_{0}\\
&1 & 0& \dots & 0& 0\\
&0 & 1& \dots & 0& 0\\
&\vdots & \vdots& \ddots & \vdots& \vdots\\
&0 & 0& \dots & 1& 0\\
\end{matrix}
\right ).
\end{equation*}
By the correspondence between polynomials and companion matrices, the eigenvalues of
$M$ are $\alpha_{1}, \dots, \alpha_{p}$, so that $M$ is diagonalizable and satisfies
\begin{align*}
&M^{m}=I,& &M^{j} \not = I,& &1 \leq j \leq m-1.&
\end{align*}
For $n \geq 1$, define
\begin{equation*}
D_{n}:=
\left (
\begin{matrix}
&a_{n-1,p-1} &a_{n-1,p-2}&
                    \dots & a_{n-1,1}&a_{n-1,0}\\
&1 & 0& \dots & 0& 0\\
&0 & 1& \dots & 0& 0\\
&\vdots & \vdots& \ddots & \vdots& \vdots\\
&0 & 0& \dots & 1& 0\\
\end{matrix}
\right ).
\end{equation*}
Thus the matrices $M$ and $D_{n}$ satisfy the conditions of Proposition \ref{tm}.
From the recurrence relation at \eqref{xrecur} we get
\begin{equation*}
\left (
\begin{matrix}
&x_{mn+i+p-1}\\
&x_{mn+i+p-2}\\
& \vdots \\
& x_{mn+i}
\end{matrix}
\right )
=\prod_{j=1}^{mn+i}D_{j}
\left (
\begin{matrix}
&x_{p-1}\\
&x_{p-2}\\
& \vdots \\
& x_{0}
\end{matrix}
\right ).
\end{equation*}
Let $F$ have the same meaning as in Proposition \ref{tm} and then
\begin{equation}\label{recurlim}
\lim_{n \to \infty}
\left (
\begin{matrix}
&x_{mn+i+p-1}\\
&x_{mn+i+p-2}\\
& \vdots \\
& x_{mn+i}
\end{matrix}
\right )
=F\, M^{i}
\left (
\begin{matrix}
&x_{p-1}\\
&x_{p-2}\\
& \vdots \\
& x_{0}
\end{matrix}
\right ).
\end{equation}

(\ref{limrecur}) now follows immediately by letting $n\to\infty$ in (\ref{pereq}).
This completes the proof.
\end{proof}

When a specific $M$ is known, (\ref{recurlim}) can sometimes be used to obtain
further relations between the different limits. This is illustrated in the following corollary.

\begin{corollary}
Let $u$ and $v$ be complex numbers, $(u,v)\not=(0,0)$ and let $\{ a_{n} \}_{n \geq 1}$,
$\{ b_{n} \}_{n \geq 1}$
 be  sequences of complex numbers such that
\begin{align*}
&\sum_{n=1}^{\infty} |a_{n}|<\infty,& &\sum_{n=1}^{\infty} |b_{n}|<\infty.&
\end{align*}
Let $\omega_{1}$ and $\omega_{2}$
be distinct roots of unity and let $m$ be the least positive integer such that
$\omega_{1}^{m}=\omega_{2}^{m}=1$.
Let the sequence $\{x_{n}\}_{n \geq 0}$ be defined by $x_{0}=u$,  $x_{1}=v$ and,
 for $n \geq 2$,
\begin{equation}\label{xrecur}
x_{n}=(\omega_{1}+\omega_{2}+a_{n-1})x_{n-1}-(\omega_{1}\omega_{2}+b_{n})x_{n-2}.
\end{equation}
Then, \\
(i) For fixed  integer $j$ the sequence $\{x_{mn+j}\}_{n \geq 0}$ is convergent;\\
(ii) If we set   $l_{j}:=\lim_{n \to \infty} x_{ mn+j}$,  then for integer $j$,
\[
l_{j+1}=(\omega_{1}+\omega_{2})\,l_{j} - \omega_{1}\omega_{2}l_{j-1};
\]
(iii) If $m$ is even and $\omega_{1}$ and $\omega_{2}$ are primitive $m$-th roots of unity, then
\[
l_{m/2+j}=-l_{j}, \hspace{30pt} 0 \leq j \leq m/2-1;
\]
(iv)For $j \in \{1,2,\dots m-2\}$, at most one of $l_{j-1}$, $l_{j}$ and $l_{j+1}$ is zero.
\end{corollary}
\begin{proof}
Define
{\allowdisplaybreaks
\begin{align*}
&M=
 \left (
\begin{matrix}
& \omega_{1}+\omega_{2}&-\omega_{1}\omega_{2} \\
&\phantom{}\\
&1&0
\end{matrix}
\right ),
\end{align*}
}
and, for $n \geq 1$, set
{\allowdisplaybreaks
\begin{align*}
&D_{n}=
\left (
\begin{matrix}
& \omega_{1}+\omega_{2}+a_{n}&-\omega_{1}\omega_{2}-b_{n} \\
&\phantom{}\\
&1&0
\end{matrix}
\right ),&
\end{align*}
}
Statement (i) follows from the $p=2$ case of Theorem \ref{t2}, since the equation
\[
t^2-(\omega_{1}+\omega_{2})t+\omega_{1}\omega_{2}=0
\]
has roots $\omega_{1}$ and $\omega_{2}$.
Statement (ii) follows immediately from Theorem \ref{t2}.
Statement (iii) follows from the fact that under the given conditions,
$
M^{m/2}=-I$ and (\ref{recurlim}) gives 
\[
\left (
\begin{matrix}
&l_{m/2+j+1}  \\
& \phantom{as} \\
&l_{m/2+j}
\end{matrix}
\right )
=
F\, M^{m/2+j}
\left (
\begin{matrix}
& u \\
&\phantom{}\\
&v
\end{matrix}
\right )
=
-F\, M^{j}
\left (
\begin{matrix}
& u \\
&\phantom{}\\
&v
\end{matrix}
\right )
=-
\left (
\begin{matrix}
&l_{j+1}  \\
& \phantom{as} \\
&l_{j}
\end{matrix}
\right ).
\]
If any two of $l_{j-1}$, $l_{j}$ and $l_{j+1}$ were zero, (iii) would then give that the third
would also be zero.  Thus
\[
\left (
\begin{matrix}
&l_{j+1}  \\
& \phantom{as} \\
&l_{j}
\end{matrix}
\right )
=
\left (
\begin{matrix}
&0  \\
& \phantom{as} \\
&0\end{matrix}
\right )
=
F\, M^{j}
\left (
\begin{matrix}
& u \\
&\phantom{}\\
&v
\end{matrix}
\right ),
\]
which is a contradiction, since $det\, F= det M = 1$ and $(u,v) \not = (0,0)$.
\end{proof}

Theorem 3.1 from \cite{ABSYZ02} follows from the above corollary, upon setting
$\omega_{1}=\exp(2 \pi i/6)$,  $\omega_{2}$ $=$ $\exp(-2 \pi i/6)$ and $b_{n}=0$ for
$n \geq 1$.

\section{A Generalization of the Ramanujan Continued Fraction}
We now study a generalization of Ramanujan's continued
fraction (\ref{R3}). As above, let $m$ be any arbitrary integer greater
than $2$, let $\omega$ be a primitive $m$-th root of unity and,
for ease of notation,
  let $\bar{ \omega} = 1/\omega$.  Define
{\allowdisplaybreaks
\begin{equation}\label{ramgen}
G(q):=
\frac{1}{1}
\-
\frac{1}{\omega + \bar{ \omega}+q}
\-
\frac{1}{\omega + \bar{ \omega}+q^2}
\-
\frac{1}{\omega + \bar{ \omega}+q^3}
\+
\cds.
\end{equation}
}

We let $P_{N}(q)/Q_{N}(q)$ denote the $N$-th approximant of $G(q)$.
From Theorem \ref{T1}, the sequence of approximants in each of the $m$ arithmetic progressions
modulo $m$ converges (set $g(x):=0$ and $f(x)=x$ in this theorem). We proceed initially along
the same path as that followed by the authors in \cite{ABSYZ02}.
We recall the $q$-binomial theorem \cite{A76}, pp. 35--36.
\begin{lemma}\label{qbin}
If
$\left [
\begin{matrix}
n\\
m
\end{matrix}
\right ]
$
denotes the Gaussian polynomial defined by
\[
\left [
\begin{matrix}
n\\
m
\end{matrix}
\right ]
:=
\left [
\begin{matrix}
n\\
m
\end{matrix}
\right ]_{q}
:=
\begin{cases}
\displaystyle{
\frac{(q;q)_{n}}{(q;q)_{m}(q;q)_{n-m}}}, &\text{ if } 0 \leq m \leq n,\\
0, &\text{ otherwise },
\end{cases}
\]
then
\begin{align}\label{qbineq}
&(z;q)_{N}= \sum_{j=0}^{N}\left [
\begin{matrix}
N\\
j
\end{matrix}
\right ]
(-1)^{j}z^{j}q^{j(j-1)/2},\\
&\frac{1}{(z;q)_{N}}
=\sum_{j=0}^{\infty}
\left [
\begin{matrix}
N+j-1\\
j
\end{matrix}
\right ]
z^{j}. \notag
\end{align}
\end{lemma}

\begin{lemma}\label{pnlem}
{\allowdisplaybreaks
\begin{equation*}
P_{N}(q) = \sum_{\substack{j,r,s \geq \, 0\\
                      r+j+s+1=N}}q^{j(j+1)/2}\omega^{r-s}
\left [
\begin{matrix}
j+r\\
j
\end{matrix}
\right ]
\left [
\begin{matrix}
j+s\\
j
\end{matrix}
\right ].
\end{equation*}
}
\end{lemma}
\begin{proof}
For $N \geq 2$, the sequence $\{P_{N}(q)\}$ satisfies
{\allowdisplaybreaks
\begin{equation}\label{precur}
P_{N}(q)=(\omega + \bar{ \omega}+q^{N-1})P_{N-1}(q)-P_{N-2}(q),
\end{equation}
}
with $P_{1}(q)=1$ and $P_{0}(q)=0$. Define
{\allowdisplaybreaks
\begin{equation*}
F(t):= \sum_{N=1}^{\infty}P_{N}(q)t^{N}.
\end{equation*}
}
If the recurrence relation \eqref{precur} is multiplied by $t^{N}$ and
summed for $N \geq 2$, we have
{\allowdisplaybreaks
\begin{equation*}
F(t)-t=(\omega + \bar{ \omega}) t F(t)+tF(tq)-t^2F(t).
\end{equation*}
}
Thus
{\allowdisplaybreaks
\begin{align*}
F(t)&=\frac{t}{(1- \omega t)(1-\bar{ \omega}t)}+
\frac{t}{(1- \omega t)(1-\bar{ \omega}t)}F(tq).
\end{align*}
}
Iterating this equation and noting that $F(0)=0$, we have that
{\allowdisplaybreaks
\begin{align*}
F(t)&=\sum_{j=0}^{\infty}
\frac{t^{j+1}q^{j(j+1)/2}}
{(\omega t; q)_{j+1}(\bar{\omega }t; q)_{j+1}}\\
&= \sum_{j,r,s=0}^{\infty}t^{j+r+s+1}q^{j(j+1)/2}\omega^{r-s}
\left [
\begin{matrix}
j+r\\
j
\end{matrix}
\right ]
\left [
\begin{matrix}
j+s\\
j
\end{matrix}
\right ],
\end{align*}
} where the last equation follows from the second equation of
(\eqref{qbineq}). The result follows upon   comparing coefficients
of $t^{N}$.
\end{proof}

\begin{lemma}\label{qnlem}
{\allowdisplaybreaks
\begin{multline*}
Q_{N}(q) = \sum_{\substack{j,r,s \geq \, 0\\
                r+j+s+1=N}}q^{j(j+1)/2}\omega^{r-s}
\left [
\begin{matrix}
j+r\\
j
\end{matrix}
\right ]
\left [
\begin{matrix}
j+s\\
j
\end{matrix}
\right ]-\\
\sum_{\substack{j,r,s \geq \, 0\\
                      r+j+s+2=N}}q^{j(j+3)/2}\omega^{r-s}
\left [
\begin{matrix}
j+r\\
j
\end{matrix}
\right ]
\left [
\begin{matrix}
j+s\\
j
\end{matrix}
\right ].
\end{multline*}
}
\end{lemma}

\begin{proof}
The proof is similar to that of Lemma \ref{pnlem}.
For $N \geq 2$, the sequence $\{Q_{N}(q)\}$
satisfies
{\allowdisplaybreaks
\begin{equation}\label{qrecur}
Q_{N}(q)=(\omega + \bar{ \omega}+q^{N-1})Q_{N-1}(q)-Q_{N-2}(q),
\end{equation}
}
with $Q_{1}(q)=Q_{0}(q)=1$. Define
{\allowdisplaybreaks
\begin{equation*}
G(t):= \sum_{N=1}^{\infty}Q_{N}(q)t^{N}.
\end{equation*}
}
If the recurrence relation \eqref{precur} is multiplied by $t^{N}$ and
summed for $N \geq 2$, we have
{\allowdisplaybreaks
\begin{equation*}
G(t)-t=(\omega + \bar{ \omega}) t G(t)+tG(tq)-t^2G(t)-t^2.
\end{equation*}
}
Thus
{\allowdisplaybreaks
\begin{align*}
G(t)&=\frac{t(1-t)}{(1- \omega t)(1-\bar{ \omega}t)}+
\frac{t}{(1- \omega t)(1-\bar{ \omega}t)}G(tq).
\end{align*}
}
Iterating this equation and noting that $G(0)=0$, we have that
{\allowdisplaybreaks
\begin{align*}
G(t)&=\sum_{j=0}^{\infty}
\frac{t^{j+1}(1-t q^{j})q^{j(j+1)/2}}
{(\omega t; q)_{j+1}(\bar{\omega }t; q)_{j+1}}\\
&= \sum_{j,r,s=0}^{\infty}t^{j+r+s+1}(1-t q^{j})q^{j(j+1)/2}\omega^{r-s}
\left [
\begin{matrix}
j+r\\
j
\end{matrix}
\right ]
\left [
\begin{matrix}
j+s\\
j
\end{matrix}
\right ],
\end{align*}
} where the last equation follows from the second equation of
(\eqref{qbineq}). The result follows upon   comparing coefficients
of $t^{N}$.
\end{proof}

\begin{lemma}\label{lemsq}
Let $\alpha$ be a primitive $m$-th root of unity and let $|q|<1$.

(i) For $j \geq 0$, $k \geq 0$ and $ i \in \mathbb{Z}$ define
{\allowdisplaybreaks
\begin{equation*}
G_{k}(\alpha, i,j,q) :=\sum_{u=0}^{mk+i}
\alpha^{u}(q^{u+1};q)_{j}.
\end{equation*}
} Then $G(\alpha, i,j,q):=\lim_{k \to \infty}G_{k}(\alpha, i,j,q)$
exists and is finite. Further, {\allowdisplaybreaks
\begin{equation}\label{Giter}
G(\alpha, i+1,j,q)=G(\alpha, i,j,q)+ \alpha^{i+1}.
\end{equation}
}

(ii)For $j \geq 0$, $k \geq 0$ and $ i \in \mathbb{Z}$ define
%{\allowdisplaybreaks
\begin{equation*}
F_{k}(\alpha, i,j,q) := \sum_{u=0}^{
\lfloor \frac{mk+i}{2} \rfloor '}
(q^{u+1};q)_{j}(\alpha^{2u-i} + \alpha^{i-2u}),
\end{equation*}
%}
where the summation
$\sum _{u=0}^{ \lfloor \frac{mk+i}{2} \rfloor '}$ means that if
$mk+i$ is even, the final term in the sum
%, corresponding to
%$u=\lfloor (mk+i)/2 \rfloor $,
is $(q^{(mk+i)/2+1};q)_{j}$, rather
 than
$2(q^{(mk+i)/2+1};q)_{j}$.

Then $F(\alpha, i,j,q):=\lim_{k \to \infty}F_{k}(\alpha, i,j,q)$
exists and is finite. Moreover,
\begin{equation}\label{FFkest} |F(\alpha,i,j,q)-F_k(\alpha,i,j,q)|\le
\frac{m2^j|q|^{{(mk+i)}/2}}{1-|q|^{{m}/{2}}.} \end{equation}
Finally,
\begin{equation}\label{Fest}
|F(\alpha,i,j,q)|\le 2^j\left(|i|+2+\frac{m|q|^{{i}/2}}{1-|q|^{{m}/{2}}}\right).
\end{equation}
\end{lemma}
Remark: As usual, we have defined the empty sum to be equal to
zero.

\begin{proof}
(i) By expanding the product $(q^{u+1})_{j}$ it be easily seen that
 \[
|(q^{u+1})_{j}-1| \leq 2^{j}|q^{u+1}|.
\]
 Since
$\sum_{u=mk+i+1}^{m(k+1)+i}\alpha^{u}=0$,
it follows that
{\allowdisplaybreaks
\begin{align*}
|G_{k+1}(\alpha, i,j,q)-G_{k}(\alpha, i,j,q)| &= \left|\sum_{u=mk+i+1}^{m(k+1)+i}
(\alpha^{u}(q^{u+1})_{j}
-\alpha^{u})\right|\\
&\leq m2^{j}|q^{mk+i}|.
\end{align*}
} This is sufficient to show that the sequence $\{G_{k}(\alpha,
i,j,q)\}_{k=0}^{\infty}$ is a Cauchy sequence and hence has a
finite limit. (\ref{Giter}) follows from the fact that
\[
G_{k}(\alpha, i+1,j,q)
=G_{k}(\alpha, i,j,q) + \alpha^{mk+i+1}(q^{mk+i+2})_{j}.
\]
Remark: The proof of (i) is not necessary for the proof of our
theorem and we give it for completeness only.

(ii) We distinguish the cases where $m$ is odd or even. If $m$ is
even then  the terms in the sequence $\{mk+i\}_{k=0}^{\infty}$ all
have the same parity, depending on the parity of $i$. Suppose $i$
is odd. Then
\begin{align}\label{feq}
F_{k+1}(\alpha, i,j,q) -F_{k}(\alpha, i,j,q) =\sum_{u=\frac{mk+i+1}{2}}^{
 \frac{m(k+1)+i-1}{2}  }
(q^{u+1};q)_{j}(\alpha^{2u-i} + \alpha^{i-2u}).
\end{align}
The collection of roots of unity, $\{\alpha^{2u-i} ,
\alpha^{i-2u}\}_{u=\frac{mk+i+1}{2}}^{ \frac{m(k+1)+i-1}{2}}$
consists of exactly two copies of the collection $\{\alpha,
\alpha^{3}, \dots , \alpha^{m-1}\}$. Thus it follows, by a similar
argument to that used in part (i), that subtracting twice the zero
sum $\sum_{v=1}^{m/2} \alpha^{2v - 1}$ from both sides of
(\ref{feq}), implies that
\[
|F_{k+1}(\alpha, i,j,q) -F_{k}(\alpha, i,j,q) | \leq m 2^{j}|q^{(mk+i+1)/2}|.
\]
This is sufficient to show that the sequence
$\{F_{k}(\alpha, i,j,q) \}$ is Cauchy and thus give the result for this case.

If both $m$ and $i$ are even, then, recalling
the definition of $F_{k}(\alpha, i,j,q)$ in the case $mk+i$ is even,
 we have that
\begin{multline}
F_{k+1}(\alpha, i,j,q) -F_{k}(\alpha, i,j,q)
=(q^{\frac{mk+i}{2}+1};q)_{j}+(q^{\frac{m(k+1)+i}{2}+1};q)_{j}\\
+\sum_{u=\frac{mk+i}{2}+1}^{
 \frac{m(k+1)+i}{2} -1 }
(q^{u+1};q)_{j}(\alpha^{2u-i} + \alpha^{i-2u}).
\end{multline}
In this case subtracting twice the zero sum
$\sum_{j=0}^{m/2-1}\alpha^{2 j}$ from both sides gives
\[
|F_{k+1}(\alpha, i,j,q) -F_{k}(\alpha, i,j,q) | \leq m
2^{j}|q|^{(mk+i)/2}.
\]
and the result follows in this case.

If $m$ is odd, the terms in the sequence $\{mk+i\}_{k=0}^{\infty}$
alternate between even and odd,
so that either
\begin{multline*}
F_{k+1}(\alpha, i,j,q) -F_{k}(\alpha, i,j,q) \\
=(q^{\frac{mk+i}{2}+1};q)_{j}
+\sum_{u=\frac{mk+i}{2}+1}^{
 \frac{m(k+1)+i-1}{2}  }
(q^{u+1};q)_{j}(\alpha^{2u-i} + \alpha^{i-2u}),
\end{multline*}
or
\begin{multline*}
F_{k+1}(\alpha, i,j,q) -F_{k}(\alpha, i,j,q) \\
=(q^{\frac{m(k+1)+i}{2}};q)_{j}
+\sum_{u=\frac{mk+i+1}{2}}^{
 \frac{m(k+1)+i}{2} -1 }
(q^{u+1};q)_{j}(\alpha^{2u-i} + \alpha^{i-2u}).
\end{multline*}
In either of these cases,  subtracting  the zero sum
$\sum_{j=0}^{m-1}\alpha^{j}$ from both sides gives
 that
\[
|F_{k+1}(\alpha, i,j,q) -F_{k}(\alpha, i,j,q) | \leq m 2^{j}|q^{(mk+i)/2}|.
\]
and the result follows once again.

Note that in all cases we have
{\allowdisplaybreaks
\begin{equation}\label{fdiffeq}
|F_{k+1}(\alpha, i,j,q) -F_{k}(\alpha, i,j,q) | \leq m 2^{j}|q^{(mk+i)/2}|.
\end{equation}
}

For brevity we let $F$ and $F_r$ denote $F(\alpha,i,j,q)$ and
$F_r(\alpha,i,j,q)$, respectively. Since $F_r\rightarrow F$, we
have for $k\ge 0$,
\[F-F_k=\sum_{r\ge k}F_{r+1}-F_r.
\]
Hence there comes \eqref{FFkest}:
\[
\begin{split}
    |F-F_k| & \le \sum_{r\ge k}|F_{r+1}-F_r|\le \sum_{r\ge k}m
2^{j}|q^{(mr+i)/2}|\\
        & \le \frac{m2^j|q|^{{(mk+i)}/2}}{1-|q|^{{m}/{2}}}.
\end{split}
\]

Now, it is easy to see directly from the definition of $F_k$ that $|F_0|\le (|i|+2)2^j$.
Putting this along with $k=0$ in the last inequality and using the triangle inequality one
more time gives (\ref{Fest}).
\end{proof}

Note that
$F(\alpha, i+m,j,q) =F(\alpha, i,j,q) $, for all $i \in \mathbb{Z}$.

\begin{lemma}\label{pqlem}
Let $1 \leq i \leq m$. With the notation of Lemma \ref{lemsq},
{\allowdisplaybreaks
\begin{equation}\label{peq}
\lim_{k \to \infty}P_{mk+i}(q)=
\sum_{j=0}^{\infty}\frac{q^{j(j+1)/2}}{(q;\,q)_{j}^{2}}F(\omega,i-j-1,j,q).
\end{equation}
}
{\allowdisplaybreaks
\begin{multline}\label{qeq}
\lim_{k \to \infty}Q_{mk+i}(q)=
\sum_{j=0}^{\infty}\frac{q^{j(j+1)/2}}{(q;\,q)_{j}^{2}}F(\omega,i-j-1,j,q)\\
\phantom{asdaaasaassasdasadssd}
-\sum_{j=0}^{\infty}\frac{q^{j(j+3)/2}}{(q;\,q)_{j}^{2}}F(\omega,i-j-2,j,q).
\end{multline}
} Moreover, the sums on the right hand side of (\ref{peq}) and
(\ref{qeq}) converge absolutely.
\end{lemma}

\begin{proof}
First note that the last assertion of the lemma follows
immediately from (\ref{Fest}).

From Lemma \ref{pnlem}, {\allowdisplaybreaks
\begin{align*}
P_{mk+i}(q)
&= \sum_{\substack{j,r,s \geq \, 0\\
                     r+j+s+1=mk+i}}q^{j(j+1)/2}\omega^{r-s}
\left [
\begin{matrix}
j+r\\
j
\end{matrix}
\right ]
\left [
\begin{matrix}
j+s\\
j
\end{matrix}
\right ]\\
&=
\sum_{\substack{j,r,s \geq \, 0\\
                     r+j+s+1=mk+i}}
\frac{q^{j(j+1)/2}}{(q)_{j}^{2}}
\omega^{r-s}
(q^{r+1})_{j}(q^{s+1})_{j}\\
&= \sum_{j=0}^{mk+i-1} \frac{q^{j(j+1)/2}}{(q)_{j}^{2}}
\sum_{r=0}^{mk+i-j-1}\omega^{2r-(mk+i-j-1)} (q^{r+1})_{j}
(q^{mk+i-j-r})_{j}\\
&=\sum_{j=0}^{mk+i-1} \frac{q^{j(j+1)/2}}{(q)_{j}^{2}}H_{k}(\omega,i,j,q),
\end{align*}
} where
{\allowdisplaybreaks
\begin{align*}
H_{k}(\omega,&i,j,q):=\sum_{r=0}^{mk+i-j-1}\omega^{2r-(mk+i-j-1)} (q^{r+1})_{j}
(q^{mk+i-j-r})_{j}\\
&=\sum_{r=0}^{ \lfloor \frac{mk+i-j-1}{2} \rfloor '}
(q^{r+1})_{j}
(q^{mk+i-j-r})_{j}(\omega^{2r-(mk+i-j-1)}+\omega^{mk+i-j-1-2r})\\
&=\sum_{r=0}^{ \lfloor \frac{mk+i-j-1}{2} \rfloor '}
(q^{r+1})_{j}
(q^{mk+i-j-r})_{j}(\omega^{2r-(i-j-1)}+\omega^{i-j-1-2r}).
\end{align*}
}
Here the summation $\sum_{r=0}^{ \lfloor \frac{mk+i-j-1}{2} \rfloor '}$
has a meaning similar to that
in Lemma \ref{lemsq}, in that if
$mk+i-j-1$ is even, then the final term is
$(q^{\frac{mk+i-j+1}{2}};q)_{j}^{2}$, rather than
$2(q^{\frac{mk+i-j+1}{2}};q)_{j}^{2}$.
The sequence $\{(q^{r+1})_{j}\}_{r=0}^{\infty}$ is bounded by $2^{j}$ and
$|\omega^{2r-(i-j-1)}+\omega^{i-j-1-2r}|\leq 2$. Thus
{\allowdisplaybreaks
\begin{align}\label{fhdiffeq}
|F_{k}(\omega, i-j&-1,j,q)-H_{k}(\omega, i,j,q)|
\leq
2^{j+1}\sum_{r=0}^{ \lfloor \frac{mk+i-j-1}{2} \rfloor }
|1-(q^{mk+i-j-r};q)_{j}|\\
&=
2^{j+1}\sum_{r=\lceil \frac{mk+i-j+1}{2} \rceil}^{ mk+i-j }
|1-(q^{r};q)_{j}|
 \leq 2^{2j+1}\sum_{r=\lceil \frac{mk+i-j+1}{2} \rceil}^{ mk+i-j }
|q|^{r} \notag \\
&\leq 2^{2j+1}\frac{|q|^{\lceil \frac{mk+i-j+1}{2} \rceil}}{1-|q|}. \notag
\end{align}
}
After applying the triangle inequality, we have that
{\allowdisplaybreaks
\begin{align*}
&\left |
P_{mk+i}
-\sum_{j=0}^{\infty}\frac{q^{j(j+1)/2}}{(q;q)_{j}^{2}}
F(\omega,i-j-1,j,q)
\right | \\
 &\leq \left |\sum_{j=mk+i}^{\infty}\frac{q^{j(j+1)/2}}{(q;q)_{j}^{2}}F(\omega,i-j-1,j,q) \right |\\
&\phantom{sadadadsada}+ \left |\sum_{j=0}^{mk+i-1}\frac{q^{j(j+1)/2}}{(q;q)_{j}^{2}}
(H_{k}(\omega, i, j, q)-F(\omega,i-j-1,j,q)) \right |\\
&\leq \sum_{j=mk+i}^{\infty}\frac{|q|^{j(j+1)/2}}{|(q;q)_{j}|^{2}}|F(\omega,i-j-1,j,q)| \\
&\phantom{sadadad}+ \sum_{j=0}^{mk+i-1}\frac{|q|^{j(j+1)/2}}{|(q;q)_{j}|^{2}}
|H_{k}(\omega, i, j, q)-F_{k}(\omega,i-j-1,j,q)| \\
&\phantom{sadadad}+ \sum_{j=0}^{mk+i-1}\frac{|q|^{j(j+1)/2}}{|(q;q)_{j}|^{2}}
|F_{k}(\omega, i-j-1, j ,q)-F(\omega,i-j-1,j,q)|.
\end{align*}
} Now apply \eqref{Fest} to the first sum, \eqref{fhdiffeq} to the second sum, and
\eqref{FFkest} to the third sum to obtain {\allowdisplaybreaks
\begin{align*}\label{fhdiffeq2}
&\left | P_{mk+i} -\sum_{j=0}^{\infty}\frac{q^{j(j+1)/2}}{(q;q)_{j}^{2}}
F(\omega,i-j-1,j,q) \right | \\
&\leq \sum_{j=mk+i}^{\infty}\frac{|q|^{j(j+1)/2}}{|(q;q)_{j}|^{2}}
2^j\left (|i-j-1|+2+\frac{m|q|^{(i-j-1)/2}}{1-|q|^{m/2}} \right ) +\\
&
%\phantom{asda}
\sum_{j=0}^{mk+i-1}\frac{|q|^{j(j+1)/2}}{|(q;q)_{j}|^{2}}
2^{2j+1}\frac{|q|^{\lceil \frac{mk+i-j+1}{2} \rceil}}{1-|q|} +
\sum_{j=0}^{mk+i-1}\frac{|q|^{j(j+1)/2}}{|(q;q)_{j}|^{2}}
\frac{m2^{j}|q|^{\frac{mk+i-j-1}{2}}}{1-|q|^{m/2}}.
\end{align*}
}
%&\leq \sum_{j=mk+i}^{\infty}\frac{|q|^{j(j+1)/2}}{|(q;q)_{j}|^{2}}
%\left (m 2^{j+1}+\frac{m\,2^{j}}{1-|q|^{m/2}} \right ) \\
%&\phantom{asda}+ \sum_{j=0}^{mk+i-1}\frac{|q|^{j(j+1)/2}}{|(q;q)_{j}|^{2}}
%2^{2j+1}\frac{|q|^{\lceil \frac{mk+i-j+1}{2} \rceil}}{1-|q|}
%+ \sum_{j=0}^{mk+i-1}\frac{|q|^{j(j+1)/2}}{|(q;q)_{j}|^{2}}
%\frac{m2^{j}|q|^{mk/2}}{1-|q|^{m/2}}.
%\end{align*}
%}
%To get the bound for the first sum, we recall that
% $F_{r}(\omega,i-j-1,j,q) = F_{r}(\omega,i',j,q)$, where $0 \leq i' \leq m-1$ and
%$i' \equiv i-j-1 (\text{mod } m )$,
%write
%%{\allowdisplaybreaks
%\begin{equation*}
%|F_{r}(\omega,i',j,q) |
%\leq
%|F_{r}(\omega,i',j,q)-F_{0}(\omega,i',j,q)| +|F_{0}(\omega,i',j,q)|,
%\end{equation*}
%}
% apply the triangle inequality and the inequality
%\eqref{fdiffeq} to $|F_{r}(\omega, i',j,q)-F_{0}(\omega,
%i',j,q)|$, let $r \to \infty$ and give a bound for $|F_{0}(\omega,
%i',j,q)|$. The  bound for the second sum comes from the inequality
%\eqref{fhdiffeq} and the third bound comes from applying the
%triangle inequality and the inequality \eqref{fdiffeq} to
%$|F_{r+k}(\omega, i-j-1,j,q)-F_{k}(\omega, i-j-1,j,q)|$ and
%letting $r \to \infty$.

The first sum is the tail of a convergent series and thus tends to
$0$ as $k \to \infty$. The third sum  is majorized by the convergent
series
{\allowdisplaybreaks
\begin{equation*}
\frac{m|q|^{(mk+i-1)/2}}{1-|q|^{m/2}}
\sum_{j=0}^{\infty}\frac{|q|^{j^{2}/2}2^{j}}{|(q;q)_{j}|^{2}},
\end{equation*}
}
and clearly  also tends to $0$ also as  $k \to \infty$. For the second sum,
{\allowdisplaybreaks
\begin{align*}
 &\sum_{j=0}^{mk+i-1}\frac{|q|^{j(j+1)/2}}{|(q;q)_{j}|^{2}}
2^{2j+1}\frac{|q|^{\lceil \frac{mk+i-j+1}{2} \rceil}}{1-|q|} \\
&= \sum_{j=0}^{\lfloor mk/2 \rfloor
-1}\frac{|q|^{j(j+1)/2}}{|(q;q)_{j}|^{2}}
2^{2j+1}\frac{|q|^{\lceil \frac{mk+i-j+1}{2} \rceil}}{1-|q|}\\
&\phantom{sadasdadadsadaSDSADSDSAsad} + \sum_{j= \lfloor mk/2
\rfloor}^{mk+i-1}\frac{|q|^{j(j+1)/2}}{|(q;q)_{j}|^{2}}
2^{2j+1}\frac{|q|^{\lceil \frac{mk+i-j+1}{2} \rceil}}{1-|q|} \\
&\leq \frac{  |q|^{ mk/4}}{1-|q|}
\sum_{j=0}^{\infty}\frac{|q|^{j(j+1)/2}}{|(q;q)_{j}|^{2}}2^{2j+1}
+ \sum_{j= \lfloor mk/2 \rfloor
}^{mk+i-1}\frac{|q|^{j(j+1)/2}}{|(q;q)_{j}|^{2}}
2^{2j+1}\frac{1}{1-|q|}.
\end{align*}
}
Both of these sums tend to $0$ as $k \to \infty$, proving \eqref{peq}. The proof
of \eqref{qeq} is virtually identical and so is omitted.
\end{proof}

We now prove Theorem \ref{ramgentheor}.
\begin{proof}[Proof of Theorem \ref{ramgentheor}]
Theorem \ref{T1} establishes the rank of the continued fraction.  The
rest of the theorem follows immediately from (\ref{ramgen}) and Lemma
\ref{pqlem}.
\end{proof}

\section{Constructing Analytic Continued Fractions with $n$ limits,
using Daniel Bernoulli's Continued Fraction}

In 1775, Daniel Bernoulli \cite{B75} proved the following result
(see, for example, \cite{K63}, pp. 11--12).
\begin{proposition}\label{pber}
Let $\{K_{0},K_{1}, K_{2},\ldots\}$ be a
sequence of complex numbers
such that  $K_{i}\not = K_{i-1}$, for $i=1,2,\ldots$.
Then  $\{K_{0},K_{1}, K_{2},\ldots\}$
is the sequence of approximants of the
continued fraction
\begin{multline}\label{ber1}
K_{0}+\frac{K_{1}-K_{0}}{1}
\+
\frac{K_{1}-K_{2}}{K_{2}-K_{0}}
\+
\frac{(K_{1}-K_{0})(K_{2}-K_{3})}
            {K_{3}-K_{1}}
\+\\
\ldots
\+
\frac{(K_{n-2}-K_{n-3})(K_{n-1}-K_{n})}
            {K_{n}-K_{n-2}}
\+
\ldots
.
\end{multline}
\end{proposition}
Trivially, if $\lim_{k \to \infty} K_{m k+ i}= L_{i}$,  for $0 \leq i \leq m-1$, where each
$L_{i}$ is different, one has a continued fraction where the approximants in each of the
$m$ arithmetic progressions modulo $m$ tend to a different limit. One easy way to use
Bernoulli's continued fraction to construct continued fractions with arbitrarily many
limits is as follows. Let  $m$ be a positive integer, $m \geq 2$. Let
$\{a_{n}\}_{n=1}^{\infty}$,  $\{c_{n}\}_{n=1}^{\infty}$,  $\{d_{n}\}_{n=1}^{\infty}$ and
 $\{e_{n}\}_{n=1}^{\infty}$
be convergent sequences with non-zero limits $a$, $c$, $d$ and $e$ respectively.
Define
\[
K_{(n-1)m+j}:= \frac{d_{n}+j\,e_{n}}{a_{n}+j\,c_{n}}
\]
for $n \geq 1$ and $0 \leq j \leq m-1$. Provided $a + j\, c \not =
0$, for  $0 \leq j \leq m-1$ and no two  consecutive terms in the
sequence $\{K_{i}\}$ are equal, then the continued fraction in
\eqref{ber1} has the sequence of  approximants
{\allowdisplaybreaks
\begin{multline*}
\biggl\{\frac{d_{1}}{a_{1}},\frac{d_{1} +e_{1}}{a_{1}+c_{1}},\frac{d_{1} +2e_{1}}
{a_{1}+2c_{1}},\cdots,\frac{d_{1} +(m-1)e_{1}}{a_{1}+(m-1)c_{1}},\cdots,\\
\frac{d_{n}}{a_{n}},\frac{d_{n} +e_{n}}
{a_{n}+c_{n}},\frac{d_{n} +2e_{n}}{a_{n}+2c_{n}},\cdots,\frac{d_{n} +(m-1)e_{n}}
{a_{n}+(m-1)c_{n}},\cdots,\biggr \}.
\end{multline*}}
Thus this continued fraction has exactly the following $m$ limits :
{\allowdisplaybreaks
\begin{equation*}
\biggl\{\frac{d}{a},\frac{d +e}{a+c},\frac{d +2e}{a+2c},\cdots,\frac{d +(m-1)e}{a+(m-1)c}\biggr\}.
\end{equation*}}
 In \cite{ABSYZ02}, the authors defined a general class of analytic continued fractions
with three limit points as follows. Let $F$ and $G$ be meromorphic functions
defined on the unit disc, $U:=\{z \in \mathbb{C}:|z|<1\}$, are analytic at the origin, and
satisfy the functional equation,
\begin{equation}\label{fgcond}
F(z)+G(z)+zF(z)G(z)=1.
\end{equation}
Further assume that $z^n, n \geq 1$, is not a pole of either $F$ or $G$.
The following theorem was proved in  \cite{ABSYZ02}.
\begin{theorem}\label{Zah}
Let $F$ and $G$ be meromorphic functions
defined on  $U$, as given above, which are analytic at the origin and satisfy
the condition \eqref{fgcond}. Then the continued fraction
\small{
\begin{multline}\label{3limcf}
\frac{1}{1}
\-
\frac{1}{1+zF(z)}
\-
\frac{1}{1+zG(z)}
\-
\frac{1}{F(z)+G(z^2)}
\-
\frac{1}{1+z^2F(z^2)}
\-
\frac{1}{1+z^2G(z^2)}
\-\\
\frac{1}{F(z^2)+G(z^3)}
\-
\cds
\-
\frac{1}{1+z^nF(z^n)}
\-
\frac{1}{1+z^nG(z^n)}
\-
\frac{1}{F(z^n)+G(z^{n+1})}
\-
\cds
\end{multline}
}
has exactly three limit points, $L_{0}$, $L_{1}$ and $L_{2}$. Moreover,
\begin{align}
L_{0}&=\frac{z}{2z-zG(z)-1},\\
L_{1}&=\frac{z + z\,G(0)-1}
   { \left(  z + z\,G(0) -1\right) \,
      \left( 1 - G(z) \right)
+\left(  z -1 \right) \,G(0) },\\
L_{2}&= \frac{1 - z\,G(0)}
   {\left( 1 - z\,G(0) \right) \,
      \left( 1 - G(z) \right) +\left(  z - 1 \right) \,\left( 1 - G(0) \right)  }.
\end{align}
\end{theorem}
We give an alternative proof, based on Proposition \ref{pber}.
\begin{proof}
If we substitute for $F$ from \eqref{fgcond} and simplify the
continued fraction, we get that the continued fraction in
\eqref{3limcf} is
 equivalent to
{\allowdisplaybreaks
\begin{multline}\label{n3limcf}
\frac{1}{1}
\-
\frac{1+z G_{1}}{1+z}
\-
\frac{1}{1}
\-
\frac{1}{1-G_{1}+(1+z G_{1})G_{2}}
\-
\frac{(1+z G_{1})(1+z^{2}G_{2})}{1+z^{2}}
\-
\cds
 \\
\-
\frac{1}{1}
\-
\frac{1}{1-G_{n-1}+(1+z^{n-1} G_{n-1})G_{n}}
\-
\frac{(1+z^{n-1} G_{n-1})(1+z^{n}G_{n})}{1+z^{n}}
\cds .
\end{multline}}
Here we use the notation $G_{n}$ for $G(z^{n})$. Define
\begin{align*}
a_{n}&:=1,\\
c_{n}&:=\frac{
\begin{matrix}
-2 + 3\,z - \left( -1 + 2\,z \right) \,{G_1} -
     \left( 1 - 2\,z + z\,{G_1} \right) \,{G_n}  \phantom{sdaASss}\\
  \phantom{sdassadaSDSAAdaasdsadsad}  + z^n\,\left( -1 + {G_1} \right) \,
      \left( z + {G_n} \right)
\end{matrix}
 }{2\,
     \left( 1 - z + \left( 1 - 2\,z + z^n \right) \,
        {G_n} + {G_1}\,\left( -1 + z +
          \left( z - z^n \right) \,{G_n} \right)  \right)
     }
,\\
d_{n}&:=\frac{1 - z + \left( -z + z^n \right) \,{G_n}}
   {1 - z + \left( 1 - 2\,z + z^n \right) \,{G_n} +
     {G_1}\,\left( -1 + z +
        \left( z - z^n \right) \,{G_n} \right) },\\
e_{n}&:=\frac{-1 + 2\,z - z^{1 + n} +
     \left( z - z^n \right) \,{G_n}}{2\,
     \left( 1 - z + \left( 1 - 2\,z + z^n \right) \,
        {G_n} + {G_1}\,\left( -1 + z +
          \left( z - z^n \right) \,{G_n} \right)  \right)
     }.
\end{align*}
In \eqref{ber1}, set $K_{0}=0$, and for $n \geq 0$, define
\begin{align}\label{berks}
K_{3n+1} &:=\frac{d_{n+1}}{a_{n+1}},\\
K_{3n+2} &:=\frac{d_{n+1}+e_{n+1}}{a_{n+1}+c_{n+1}}, \notag \\
K_{3n+3} &:=\frac{d_{n+1}+2e_{n+1}}{a_{n+1}+2c_{n+1}}. \notag
\end{align}
Then the continued fraction in  \eqref{ber1}
 simplifies to give \eqref{n3limcf}.
If we simplify the expressions on the right side of \eqref{berks},
we have that this continued fraction
has exactly the sequence of approximants
\begin{multline}
\bigg \{ \frac{1 - z + \left( -z + z^n \right) \,{G_n}}
   {1 - z + \left( 1 - 2\,z + z^n \right) \,{G_n} +
     {G_1}\,\left( -1 + z +
        \left( z - z^n \right) \,{G_n} \right) },\\
  \frac{1 - z^{1 + n} + \left( -z + z^n \right) \,{G_n}}
   {-\left( z\,\left( -1 + z^n \right)  \right)  +
     \left( 1 - 2\,z + z^n \right) \,{G_n} +
     {G_1}\,\left( -1 + z^{1 + n} +
        \left( z - z^n \right) \,{G_n} \right) },\\
  - \frac{z\,\left( -1 + z^n \right) }
     {-1 + 2\,z - z^{1 + n} +
       z\,\left( -1 + z^n \right) \,{G_1}}
 \bigg \}_{n=1}^{\infty} .
\end{multline}
Finally, we let $n \to \infty$ to get the three limits, noting that
$z^n \to 0$ and $G_{n}=G(z^n) \to G(0)$.
\end{proof}
It is not difficult to create analytic continued fractions with $m$ limit points,
where $m$ is an integer, $m \geq 2$.  For example,  Theorem \ref{t4} provides
a continued fraction which is less convoluted in its contrivance than
the one in Theorem \ref{Zah}.  It is clear that Bernoulli's continued fraction
can be used to construct many  similar examples.

\begin{proof}[Proof of Theorem \ref{t4}]
In  \eqref{ber1}, put, for $i \geq 0$,
\[
K_{i}=G(z^{i+1})- G(z) + i (\text{mod } m). \]
 The sequence $\{K_{i}\}_{i=0}^{\infty}$ has
exactly the $m$ limits stated in the theorem. For $i \geq 1$,
{\allowdisplaybreaks
\begin{equation}
K_{i}-K_{i-1}=
\begin{cases}
G(z^{i+1}) -G(z^{i})-m+1, \text{ for } i \equiv 0 (\text{mod } m ),\\
\phantom{as}\\
G(z^{i+1}) -G(z^{i})+1,  \text{ otherwise.}
\end{cases}
\end{equation}
} 
\eqref{gcon} gives that $K_{i}-K_{i-1} \not = 0$ for $i \geq 1$. The
continued fraction \eqref{bernoucf} above is simply Bernoulli's continued fraction
\eqref{ber1} for the stated sequence $\{K_{i}\}_{i=0}^{\infty}$.
\end{proof}

\section{Concluding Remarks}
It should be noted that the condition in Theorem \ref{T1} that $\omega_{1}$ and
$\omega_{2}$ be distinct is necessary, since otherwise the matrix $M$ cannot be diagonalized.
Moreover substituting $\omega_{1}=\omega_{2}=1$ (so $m=1$) into the continued fraction $G$ in
Theorem \ref{T1} does not necessarily give that $G$ has $m=1$ limit. An example of
how dropping the condition that $\omega_{1}$ and
$\omega_{2}$ be distinct can lead to a false result is provided by the continued fraction
\[
K_{n=1}^{\infty} \frac{-1 - 4/(4n^2-1)}{2} =2K_{n=1}^{\infty} \frac{-1/4 - 1/(4n^2-1)}{1}.
\]
Our Theorem \ref{T1} without the condition that $\omega_{1}$ and
$\omega_{2}$ be distinct would predict that the first continued fraction above, and hence the
second, would have one limit ($m=1$) and hence converge, but the second continued
fraction diverges generally (\cite{LW92}, page 158).

We have not encountered the functions $G(\alpha, i,j,q)$ and
$F(\alpha, i,j,q)$ from Lemma \ref{lemsq} elsewhere. It would be
interesting to know something of their properties. Do they have
product representations? Since $F(\alpha, i,j,q)=F(\alpha,
i+m,j,q)$ for all $i \in \mathbb{Z}$, is there a simple expression
relating $F(\alpha, i,j,q)$ and $F(\alpha, i+1,j,q)$, similar to
the expression in \eqref{Giter} relating $G(\alpha, i,j,q)$ and
$G(\alpha, i+1,j,q)$?
 Is there a formula relating $F(\alpha, i,j,q)$
and $F(\alpha, i,j+1,q)$?

Our formulas for the $m$ limits in Theorem \ref{ramgentheor} lack
the simplicity of Ramanujan's for the continued fraction with
three limit points. Do the quotients of series on the left side of
(\ref{ramgentheoreq}) have expressions in terms of infinite
products?

 \allowdisplaybreaks{

}
\end{document}